# Fault-tolerant Identifying Codes in Special Classes of Graphs

Devin C. Jean

Computer Science Department

Vanderbilt University

`devin.c.jean@vanderbilt.edu`

Suk J. Seo

Computer Science Department

Middle Tennessee State University

`Suk.Seo@mtsu.edu`

**Abstract**

A detection system, modeled in a graph, is composed of "detectors" positioned at a subset of vertices in order to uniquely locate an "intruder" at any vertex. *Identifying codes* use detectors that can sense the presence or absence of an intruder within distance one. We introduce a fault-tolerant identifying code called a *redundant identifying code*, which allows at most one detector to go offline or be removed without disrupting the detection system. In real-world applications, this would be a necessary feature, as it would allow for maintenance on individual components without disabling the entire system. Specifically, we prove that the problem of determining the lowest cardinality of an identifying code for an arbitrary graph is NP-complete, we determine the bounds on the lowest cardinality for special classes of graphs, including trees, cylinders, and cubic graphs.

## 1 Introduction

Let $G = (V(G), E(G))$ be a (simple) graph, with vertices $V(G)$ and edges $E(G)$, modelling a system or facility with detectors to recognize a possible problem, traditionally referred to as an "intruder". For example, the vertices of the graph can represent sections of a shopping mall, the intruder could be a shoplifter, and the detectors can be video surveillance equipment or motion, magnetic, or RFID sensors. The goal is to identify the exact location/vertex of the intruder by placing the minimum number of detectors in the facility/graph. To represent the capabilities of the sensor(s) placed at a point $v \in V(G)$, we associate each sensor, $\rho$, at location $v$ with a detection region $R_\rho(v) \subseteq V(G)$, where $\rho$ can detect the presence or absence of an intruder anywhere in $R_\rho(v)$. The vertex $v$ itself is associated with a set of detection regions, $R(v)$, which is simply the set of $R_\rho(v)$ for each sensor $\rho$ at position $v$. Note that when $|R(v)| = 1$, $R_\rho(v)$ has also been referred to as the "watching zone" of $v$ in other papers, but with different notation [1].

**Definition 1.** *Let $G$ be a graph and $v \in V(G)$. The* open neighborhood *of $v$, denoted $N(v)$, is the set of all vertices adjacent to $v$, $\{w \in V(G) : vw \in E(G)\}$.*

**Definition 2.** *Let $G$ be a graph and $v \in V(G)$. The* closed neighborhood *of $v$, denoted $N[v]$, is the set of all vertices adjacent to $v$ as well as $v$ itself, $N(v) \cup \{v\}$.*

Many types of detection systems with various properties have been explored throughout the years. One such system is the *Locating-Dominating (LD) set*, where each detector can sense the presence of an intruder in its closed neighborhood but also has the ability to distinguish the vertex itself from its neighbors; that is, $R(v) = \{\{v\}, N(v)\}$ [12]. Another type of distinguishing set that has been explored is the *open-locating-dominating (OLD) set*, which is based on LD but removed the self-distinguishing property; that is, $R(v) = \{N(v)\}$ [10]. Of particular interest in this paper are *identifying codes (ICs)*, where $R(v) = \{N[v]\}$ [8]. Over 440 papers have been published on these detection systems and other related concepts [9].

Detection systems are useful in modeling security systems and automated fault detection in networked systems; thus, it is often the case that we want some level of fault tolerance guaranteed in the system. Many different forms of fault tolerant detection systems have been explored, including the ability to correct false negative or false positive signals from the sensors. Identifying codes were introduced by Karpovsky et. al in 1998 [8]; in this paper, we will introduce *redundant identifying codes (RED:ICs)*, which allow at most one detector to go offline or be removed without disrupting the system. To the best of our knowledge, this is the first paper to consider fault-tolerant identifying codes.

Detector-based systems commonly use terminology such as "dominated" or "distinguished", whose definitions vary depending on the sensors' capabilities. The following definitions are specifically for identifying codes and their fault-tolerant variants; assume that $S \subseteq V(G)$ is the set of detectors.

**Definition 3.** *A vertex, $v \in V(G)$, is $k$-dominated if $|N[v] \cap S| = k$.*

**Definition 4.** *Two vertices, $u, v \in V(G)$, are $k$-distinguished if $|(N[u] \cap S) \triangle (N[v] \cap S)| \geq k$.*

We will also use terms such as "at least $k$-dominated" to denote $l$-dominated for some $l \geq k$.

**Theorem 1.** *A detector set, $S \subseteq V(G)$, is an IC if and only if each vertex is at least 1-dominated and all pairs are 1-distinguished.*

**Theorem 2.** *A detector set, $S \subseteq V(G)$, is a RED:IC if and only if each vertex is at least 2-dominated and all pairs are 2-distinguished.*

Theorem 2 was stated more generally by Slater [11, 13]; in this paper, it is specialized for RED:IC.

For finite graphs, we use the notations $IC(G)$, and $RED:IC(G)$ to denote the cardinality of the smallest possible such sets on graph $G$, respectively. For infinite graphs, we measure via the *density* of the subset, which is defined as the ratio of the size of the subset to the size of the whole set [7, 11]. We use the notations $IC\%(G)$, and $RED:IC\%(G)$ to denote the lowest density of any possible such set on $G$ [7, 11]. Note that density is also defined for finite graphs.

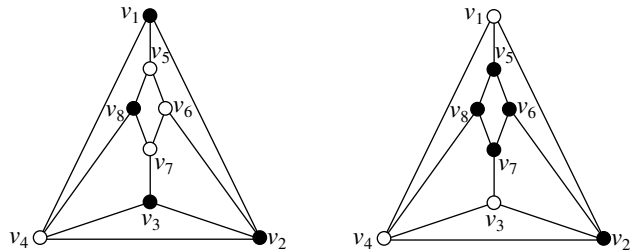

Figure 1: Optimal IC (a) and RED:IC (b) sets on $G_8$. Shaded vertices represent detectors.

Figure 1 shows an example IC and RED:IC on $G_8$. In the IC set (a), we see that $N[v_1] \cap S = \{v_1, v_2\}$, $N[v_2] \cap S = \{v_1, v_2, v_3\}$, $N[v_6] \cap S = \{v_2\}$, and so on; each set has at least one item, so every vertex is at least 1-dominated. For brevity, let $\triangle_{a,b} = (N[a] \cap S) \triangle (N[b] \cap S)$. We see that $\triangle_{v_1,v_2} = \{v_5, v_6\}$, $\triangle_{v_1,v_3} = \{v_5, v_7\}$, $\triangle_{v_5,v_7} = \{v_5, v_7\}$, and so on; each has at least one item, so all vertex pairs are 1-distinguished. Therefore, Theorem 1 yields that (a) is an IC. We can perform a similar analysis on the vertices of (b) to see that all vertices are at least 2-dominated and all pairs are 2-distinguished, so by Theorem 2 it is a RED:IC.

Thus, we see that $IC(G_8) = 4$ and $RED{:}IC(G_8) = 5$. If we would prefer to use densities, we also have that $IC\%(G_8) = \frac{1}{2}$ and $RED{:}IC\%(G_8) = \frac{5}{8}$.

In the following section, we prove that the problem of determining the minimum cardinality of $RED{:}IC(G)$ for an arbitrary graph $G$ is NP-complete. In Section 3 we show existence criteria of the redundant identifying code for general graphs, and determine the lower bound on the minimum density of RED:IC, which is a tight bound when $n$ is even. In Section 4 we explore several special classes of graph—including ladders, trees, and cubic graphs—and find lower and upper bounds of $RED{:}IC(G)$, some tight bounds, and several extremal families of graphs with minimum and maximum value. Section 5 concludes with our initial results on some infinite grids and future direction.

## 2  NP-completeness of RED:IC

It has already been shown that many graphical parameters related to detection systems, such as finding the cardinality of the smallest IC, LD, or OLD sets, are NP-complete problems [2, 3, 4, 10]. We will now show that the problem of determining the smallest RED:IC set is also NP-complete. For additional information about NP-completeness, see Garey and Johnson [6].

**3-SAT**
**INSTANCE:** Let $X$ be a set of $N$ variables. Let $\psi$ be a conjunction of $M$ clauses, where each clause is a disjunction of three literals from distinct variables of $X$.
**QUESTION:** Is there is an assignment of values to $X$ such that $\psi$ is true?

**Redundant Identifying Code (RED-IC)**
**INSTANCE:** A graph $G$ and integer $K$ with $2 \leq K \leq |V(G)|$.
**QUESTION:** Is there a RED:IC set $S$ with $|S| \leq K$? Or equivalently, is $RED{:}IC(G) \leq K$?

**Theorem 3.** *The RED-IC problem is NP-complete.*

*Proof.* Clearly, RED-IC is NP, as every possible candidate solution can be generated nondeterministically in polynomial time (specifically, $O(n)$ time), and each candidate can be verified in polynomial time using Theorem 2. To complete the proof, we will now show a reduction from 3-SAT to RED-IC.

Let $\psi$ be an instance of the 3-SAT problem with $M$ clauses on $N$ variables. We will construct a graph, $G$, as follows. For each variable $x_i$, create an instance of the $F_i$ graph (Figure 2); this includes a vertex for $x_i$ and its negation $\overline{x_i}$. For each clause $c_j$ of $\psi$, create a new instance of the $H_j$ graph (Figure 2). For each clause $c_j = \alpha \vee \beta \vee \gamma$, create an edge from the $c_j$ vertex to $\alpha$, $\beta$, and $\gamma$ from the variable graphs, each of which is either some $x_i$ or $\overline{x_i}$; for an example, see Figure 3. The resulting graph has precisely $8N + 3M$ vertices and $8N + 5M$ edges, and can be constructed in polynomial time.

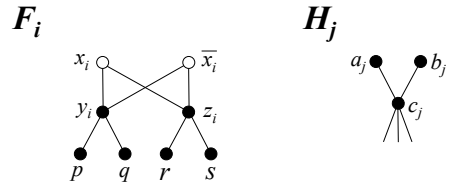

Figure 2: Variable and Clause graphs

Suppose $S \subseteq V(G)$ is an optimal (minimum) RED:IC on $G$. By Theorem 2, every vertex must be 2-dominated; thus, we require at least $6N + 3M$ detectors, as shown by the shaded vertices in Figure 2. For each $H_j$, we see that $a_j$ and $c_j$ are not distinguished unless $c_j$ is adjacent to at least one additional detector vertex. Similarly, in each $F_i$ we see that $y_i, p_i$ and $z_i, r_i$ are not distinguished unless $\{x_i, \overline{x_i}\} \cap S \neq \varnothing$. Thus, we find that $|S| \geq 7N + 3M$; if $|S| = 7N + 3M$, then $|\{x_i, \overline{x_i}\} \cap S| = 1$ and $c_j$ must be dominated by one of its three neighbors in the $F_i$ graphs, so $\psi$ is satisfiable.

Now suppose we have a solution to the 3-SAT problem. For each variable, $x_i$, if $x_i$ is true then we let the vertex $x_i \in S$; otherwise, we let $\overline{x_i} \in S$ (in addition to the $6N + 3M$ detectors that are required for 2-domination). Because $\psi$ was satisfiable by assumption, each $c_j$ must be adjacent to at least one additional detector vertex from the $F_i$ graphs. □

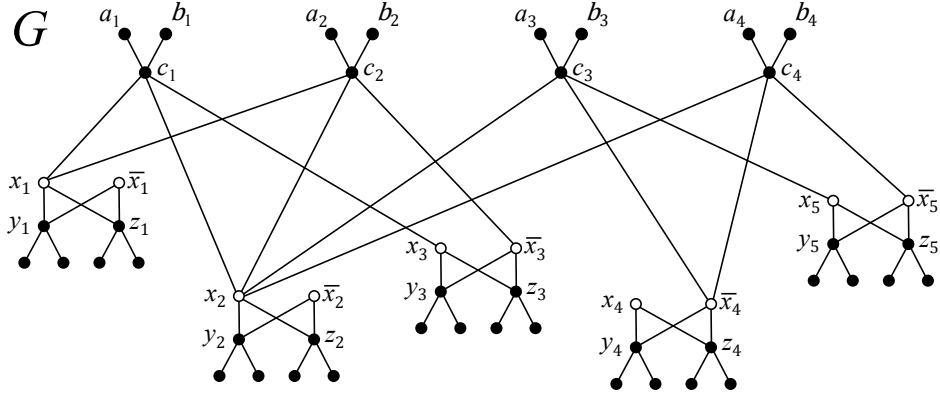

Figure 3: Construction of $G$ from $(x_1 \vee x_2 \vee x_3) \wedge (x_1 \vee x_2 \vee \overline{x_3}) \wedge (x_2 \vee \overline{x_4} \vee x_5) \wedge (x_2 \vee \overline{x_4} \vee \overline{x_5})$

## 3  Existence of RED:IC and bounds on RED:IC(G)

**Definition 5.** *If $u \in V(G)$ has $N(u) = \{v\}$, then $u$ is called a* leaf vertex *and $v$ is called a* support vertex.

**Definition 6.** *If a vertex, $v$, is neither a leaf nor support vertex, it is called a* pure interior vertex.

**Definition 7.** *[5] Two distinct vertices $u, v \in V(G)$ are said to be* twins *if $N[u] = N[v]$ (* closed twins*) or $N(u) = N(v)$ (* open twins*).*

From Theorem 2, we know that each pair of vertices must be 2-distinguished; if $u$ and $v$ are closed twins, then they cannot be distinguished, so no RED:IC exists. We also see that all support and leaf vertices must be detectors in order to 2-dominate the leaves; if a support vertex is not at least 4-dominated, then it will not be distinguished from its leaves. Thus, we arrive at the following:

**Observation 1.** *RED:IC exists only if there are no closed twins and every support vertex, $v$, has $deg(v) \geq 3$.*

**Observation 2.** *There is no graph with RED:IC(G) $\leq 3$.*

**Observation 3.** *The smallest graph with RED:IC(G) is a claw graph ($K_{1,3}$) and $C_4$ with RED:IC(G) = 4.*

**Theorem 4.** *Let $G$ be connected with $n \geq 4$. RED:IC exists if and only if there are no closed twins, every support vertex is at least degree three, and every triangle $abc \in G$ has $|N[a] \triangle N[b]| \geq 2$.*

*Proof.* Let $S = V(G)$ be a set of detectors; because $G$ is connected and $|V(G)| \geq 2$, every vertex is at least 2-dominated. We will show that each $v \in V(G)$ is distinguished from every other vertex $u \in V(G)$. If $uv \notin E(G)$, then $u$ and $v$ are distinguished by themselves; otherwise, we assume $uv \in E(G)$. By hypothesis, $u$ and $v$ are not closed twins; without loss of generality, let $w_1 \in N(u) - N[v]$. If $v$ is a leaf, then $deg(u) \geq 3$ by hypothesis, so $u$ and $v$ are distinguished by the neighbors of $u$; otherwise $\exists w_2 \in N(v)$. If $w_2 \notin N(u)$, then $u$ and $v$ are distinguished by $w_1$ and $w_2$; otherwise $uvw_2$ is a triangle, so by hypothesis $u$ and $v$ are distinguished, meaning $S$ is a RED:IC. For the converse, suppose one of the properties is not met. If $u, v \in V(G)$ are closed twins, then they cannot be distinguished, so no RED:IC exists. If $v$ is a support vertex with $deg(v) \leq 2$, then $v$ and its leaf cannot be distinguished, so no RED:IC exists. Finally, if there is a triangle $abc \in G$ with $|N[a] \triangle N[b]| \leq 1$, then $a$ and $b$ cannot be 2-distinguished, so no RED:IC exists, completing the proof. □

Based on Theorem 2, we can easily construct an algorithm to test if a connected graph $G$ with $n \geq 4$ has a RED:IC set: simply check that for any $u, v \in V(G)$, $|N[u] \triangle N[v]| \geq 2$, which can be done in $\mathcal{O}(m\Delta(G))$ time in the worst case if the graph input is an adjacency list. From Observation 1 and Theorem 4, if $G$ is triangle-free, we need only ensure that support vertices are at least degree three and that there are no closed twins; if the input is a sorted adjacency list, this can be done in $\mathcal{O}(n\Delta(G))$ time by iterating over each vertex and storing the closed neighborhoods in a set. Finally, if $G$ is a tree, we need only check that every support vertex is at least degree three, which can be done in $\mathcal{O}(n)$ time if the input is an adjacency list.

Next, we consider a lower bound on the value of RED:IC($G$). We start by analyzing the maximum size of a graph with IC($G$) = $k$. We know that every vertex must be at least 1-dominated and all pairs must be 1-distinguished (non-identical), which means that all non-empty combinations of the total $k$ detectors are valid codewords (closed neighborhoods of detectors).

**Observation 4.** *If $IC(G) \leq k$, then $|V(G)| \leq 2^k - 1$.*

**Theorem 5.** *If $RED{:}IC(G) \leq k$, then $|V(G)| \leq 2^{k-1} - 1$.*

*Proof.* Suppose we have a RED:IC, $S \subseteq V(G)$, with $|S| \leq k$. Thus, by definition, there exists an IC $S'$ with $|S'| \leq k-1$. Observation 4 gives us that $|V(G)| \leq 2^{k-1} - 1$. □

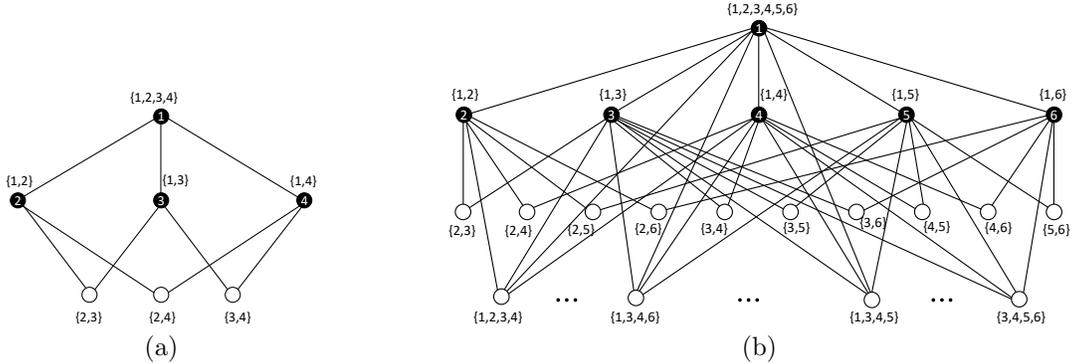

Figure 4: Extremal family of graphs with largest $n$, with examples of where RED:IC($G$) is 4 (a) and 6 (b)

For $k = 2j$, start with a star graph, $K_{1,k}$, where every vertex is a detector. We then add an additional $\binom{k}{2} - (k-1)$ non-detectors which are adjacent to a distinct pair of detectors (keeping in mind there are already $k-1$ detector vertices which have a pair including themselves), an additional $\binom{k}{4}$ non-detectors which are adjacent to distinct sets of 4 detectors, an additional $\binom{k}{6}$ non-detectors which are adjacent to distinct sets of 6 detectors, and so on through $\binom{k}{k-2}$ (as $\binom{k}{k}$ was already created at the beginning). Because only even numbers of detector neighbors were chosen, every vertex will be at least 2-dominated and 2-distinguished. Then, the total number of vertices is thus $\binom{k}{2} + \binom{k}{4} + \ldots + \binom{k}{k} = 2^{k-1} - 1$ (see Equation 1). We see that this value matches the theoretical upper bound established by Theorem 5. This construction yields an infinite family of extremal graphs with largest $n$ for any even value of RED:IC, $k$; example graphs for $k = 4$ and $k = 6$ are shown in Figure 4.

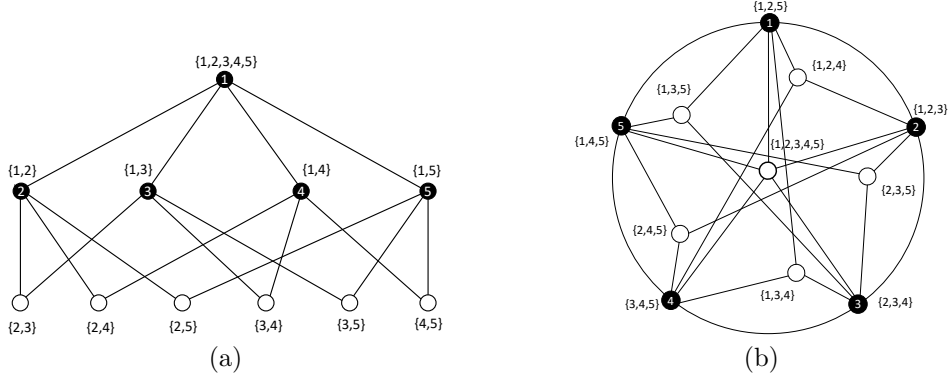

Figure 5: Two constructions of graphs with large $n$ that have RED:IC$(G) = 5$

$$2^n = (1+1)^n = \sum_{k=0}^{n} \binom{n}{k}$$
$$0 = (1-1)^n = \sum_{k=0}^{n} \binom{n}{k}(-1)^k \tag{1}$$

For $k = 2j+1$, we again start with a star graph, $K_{1,k}$, where every vertex is a detector. We add non-detectors in a similar fashion to the case when $k = 2j$, starting with $\binom{k}{2} - (k-1)$ non-detectors, but ending with $\binom{k}{k-3}$ (as non-detectors representing $\binom{k}{k-1}$ will not be distinguished from the $\binom{k}{k}$ detector). This construction is shown for $k = 5$ in Figure 5 (a). Then, the total number of vertices is $\binom{k}{2} + \binom{k}{4} + \ldots + \binom{k}{k-3} + \binom{k}{k} = 2^{k-1} - k$ (see Equation 1).

Alternatively, for $k = 2j+1$, we can start with a cycle on $k$ vertices, $C_k$, where every vertex is a detector. Add an additional ($\binom{k}{3}$ - k) non-detectors that are adjacent to distinct set of three detectors. Add $\binom{k}{5}$ non-detectors so that they are adjacent to distinct set of 5 detectors. Add $\binom{k}{7}$ non-detectors so that they are adjacent to distinct set of 7 detectors, and so on through $\binom{k}{k}$. An example of this construction for $k = 5$ is given by Figure 5 (b). Because only odd numbers of detector neighbors were chosen, every vertex will be at least 2-dominated and 2-distinguished. Then, the total number of vertices is $\binom{k}{3} + \binom{k}{5} + \ldots + \binom{k}{k-2} + \binom{k}{k} = 2^{k-1} - k$.

There is no RED:IC for a complete graph, as every vertex is a closed twin. There are $p = \lfloor \frac{n}{2} \rfloor$ distinct pairs of closed twins in $K_n$. For any of the $p$ pairs of twins, $u$ and $v$, we can remove the $uv$ edge; this makes them no longer closed twins with one another, and does not affect other vertices. If $n$ is even, removing the $p$ edges corresponding to the $p$ distinct pairs of twins results in a complete multipartite graph with $p$ parts, each of size 2. This graph must necessarily have RED:IC$(G) = n$ because every vertex is an open twin, and $G$ has the maximum number of edges that RED:IC allows, as only the necessary $p$ edges were removed.

From Theorem 5, we see that if $S$ is a RED:IC of size $k$, then $n \leq 2^{k-1} - 1$, from which we see that $\log_2(n+1) + 1 \leq k$. This gives us the following corollary:

**Corollary 1.** *If $G$ has a RED:IC, then $\lceil \log_2(n+1) \rceil + 1 \leq$ RED:IC$(G) \leq n$.*

## 4 Special classes of graphs

**Observation 5.** *If $S$ is a RED:IC set and $v \in V(G)$ is adjacent to a degree 3 support vertex, then $v \in S$.*

**Observation 6.** *If $S$ is a RED:IC set and $vuw$ is a path in $G$ where $u$ and $w$ are degree 2, then $v \in S$.*

From Observation 1, finite paths do not have RED:IC. From Observation 6, we see that the infinite path and cycles with $n \geq 4$ require RED:IC%$(G) = 1$.

### 4.1 Trees

Because a tree on $n \geq 4$ vertices is closed-twin free and triangle free, Observation 1 and Theorem 4 give us the following corollary:

**Corollary 2.** *Let $T$ be a tree with $n \geq 4$. RED:IC exists if and only if every support vertex, $v$, has $deg(v) \geq 3$.*

**Upper bound on RED:IC($T_n$) for finite trees**

By the requirement of 2-domination, a star graph $K_{1,n-1}$, has RED:IC%$(K_{1,n-1}) = 1$. In fact, from Figure 7, we see that all trees of order $4 \leq n \leq 8$ have RED:IC$(T) = n$. We will now generalize the requirements for these extremal trees.

**Theorem 6.** *If $T$ is a tree of order $n \geq 4$ with a RED:IC, then RED:IC$(T) = n$ if and only if each vertex, $v \in V(G)$, is a leaf vertex, support vertex, or satisfies Observation 5 or 6.*

*Proof.* Clearly, if $T$ satisfies the four properties in the theorem statement, then RED:IC$(T) = n$. For the converse, suppose that some $v \in V(T)$ does not satisfy any of the four properties; we will show that $V(T) - \{v\}$ is a RED:IC. By hypothesis, we know a RED:IC exists, which by Observation 1 implies that all support vertices are at least degree 3. Because $v$ is not a leaf vertex, $deg(v) \geq 2$; let $w \in N(v)$, let $T'$ be the subtree of $T - \{v\}$ containing $w$, and let $n' = |V(T')|$. We will show that $S' = V(T')$ is a RED:IC for $T'$, meaning $V(T) - \{v\}$ is a RED:IC for the original graph, $T$. Because $v$ is not a support vertex, $deg(w) \geq 2$; let $z \in N(w) - \{v\}$. If $deg(w) = deg(z) = 2$, we contradict that $v$ does not satisfy Observation 6; thus, we assume that $deg(w) \geq 3$ or $deg(z) \neq 2$. Suppose $deg(w) = 2$; then require $deg(z) \neq 2$. If $deg(z) = 1$, then $w$ is a support vertex which is not at least degree 3, a contradiction; otherwise, we assume that $deg(z) \geq 3$. We see that $n' \geq 4$ and every support vertex in the restricted graph $T'$ is at least degree 3, so Theorem 4 yields that $V(T')$ is a RED:IC on $T'$, and we are done. Otherwise, $deg(w) \geq 3$. If $deg(w) \geq 4$, we again see that $n' \geq 4$ and all support vertices in $T'$ are at least degree 3, so we are done; thus, we can assume that $deg(w) = 3$. If $w$ is a support vertex, then we contradict that $v$ does not satisfy Observation 5; thus, we assume that $w$ is not a support vertex, meaning $\forall z_i \in N(w), deg(z_i) \geq 2$. Thus, we again see that $n' \geq 4$ and all support vertices in $T'$ maintain degree at least 3, so we are done. Therefore, in any case, we find that $S' = V(T')$ is a RED:IC for $T'$, so $V(T) - \{v\}$ is a RED:IC for $T$, completing the proof. □

**Lower bound on RED:IC($T_n$) for finite trees**

From the cubic graph section, we know RED:IC$(T_\infty)$ for infinite 3-regular tree has $\frac{4}{7}$. However, for finite trees the lower bound on RED:IC(T) is much higher as shown in the next theorem.

**Theorem 7.** *If a tree, $T$, of order $n \geq 4$ has RED:IC, then $\lceil \frac{4}{5}(n+1) \rceil \leq RED:IC(T) \leq n$.*

*Proof.* Suppose $S$ is a RED:IC for $T$ and let $j = |V(G) - S|$ be the number of non-detectors; then, $n = |S| + j$. Because $T$ is acyclic and each non-detector must be at least 2-dominated, we know that there must be at least $j+1$ connected components of detectors. Because $S$ must be a RED:IC on the graph generated by $S$, each connected component of detectors must have a RED:IC, meaning the minimum size of each detector component is four. Thus, $|S| \geq 4j+4$; by rearranging terms, we see that $|S| \geq 4(n-|S|)+4$, from which we find that $|S| \geq \frac{4}{5}(n+1)$. We know that RED:IC$(T) \in \mathbb{N}$, so we can strengthen this to $|S| \geq \lceil \frac{4}{5}(n+1) \rceil$, completing the proof. □

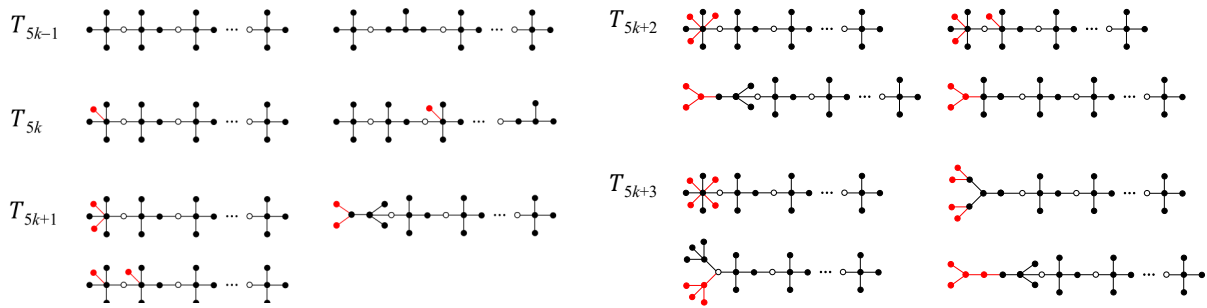

Figure 6: Extremal trees with RED:IC$(T) = \lceil \frac{4}{5}(n+1) \rceil$. Red vertices denote "excess" detectors.

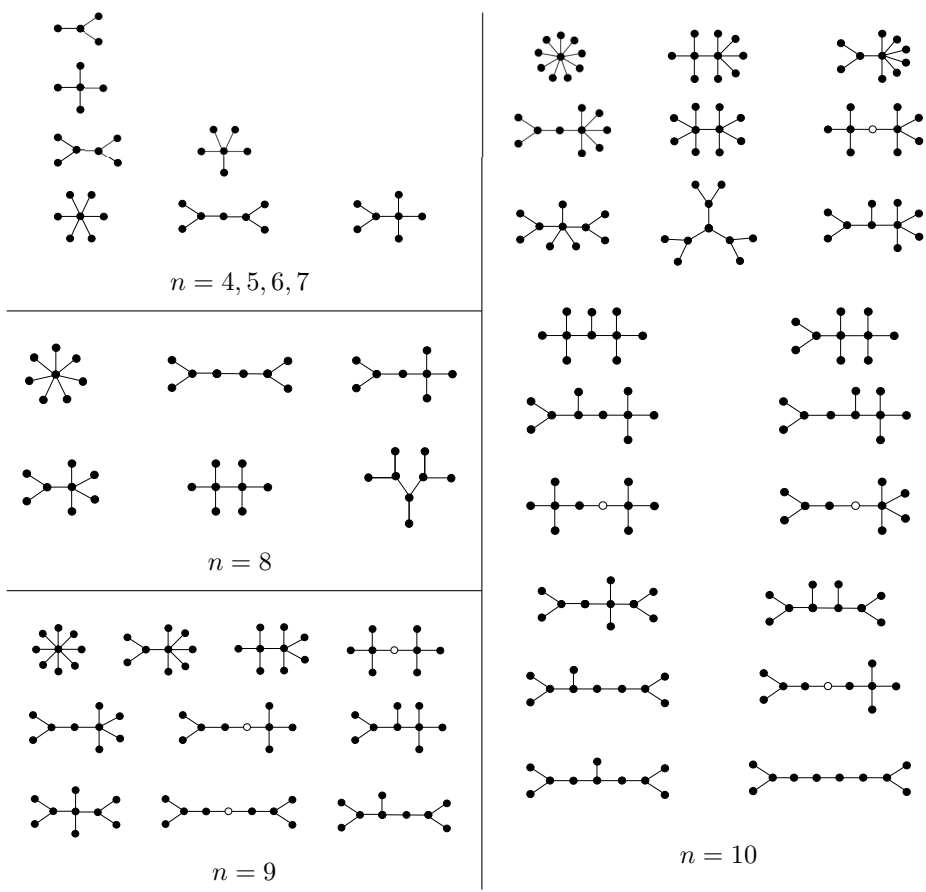

Figure 7: All trees of order $n \leq 10$ with RED:IC

Figure 6 gives examples of extremal trees on $n \geq 4$ vertices with RED:IC$(T_n) = \lceil \frac{4}{5}(n+1) \rceil$. In general, we see that when there are $j$ non-detectors, there must be $p = j+1$ detector components (or up to $p = j+2$ if $n = 5k+3$)—each detector component can be selected from Figure 7 on at most 8 vertices—and the total number of detectors must be $4p + [(n+1) \bmod 5]$; edges between non-detectors and detectors can be added arbitrarily so long as the result is a tree. Note that this is essentially the same as starting with an extremal tree on $5k-1$ vertices, as in Figure 6, and adding an additional $(n+1) \bmod 5$ "excess" detectors, so long as the placements do not cause RED:IC to no longer exist.

Figure 7 shows optimal RED:ICs on all trees of order $n \leq 10$ for which RED:IC exists. Table 1 provides a summarized view of the number of trees with a given RED:IC value for $n \leq 17$ vertices.

| $n$ | 4 | 5 | 6 | 7 | 8 | 9 | 10 | 11 | 12 | 13 | 14 | 15 | 16 | 17 |
|---|---|---|---|---|---|---|---|---|---|---|---|---|---|---|
| trees | 2 | 3 | 6 | 11 | 23 | 47 | 106 | 235 | 551 | 1301 | 3159 | 7741 | 19320 | 48629 |
| with RED:IC | 1 | 1 | 2 | 3 | 6 | 10 | 21 | 39 | 82 | 167 | 360 | 766 | 1692 | 3726 |
| RED:IC$(T_n) = n-2$ | 0 | 0 | 0 | 0 | 0 | 0 | 0 | 0 | 0 | 0 | 13 | 29 | 96 | 287 |
| RED:IC$(T_n) = n-1$ | 0 | 0 | 0 | 0 | 0 | 3 | 4 | 10 | 24 | 64 | 130 | 323 | 744 | 1731 |
| RED:IC$(T_n) = n$ | 1 | 1 | 2 | 3 | 6 | 7 | 17 | 29 | 58 | 103 | 217 | 414 | 852 | 1708 |

Table 1: Results on RED:ICs for trees

## 4.2 Ladders and cylinders

**The infinite ladder**

**Theorem 8.** *The infinite ladder graph has RED:IC%*$(P_2 \square P_\infty) = \frac{2}{3}$.

*Proof.* Figure 9 (c) and (d) give a family of cylinders with RED:IC%$(G) = \frac{2}{3}$; each of these solutions can be tiled infinitely to produce a RED:IC on $P_2 \square P_\infty$ with density $\frac{2}{3}$. To prove this is optimal, we need only show that $\frac{2}{3}$ is a lower bound for the minimum density. To proceed, we will look at an arbitrary non-detector vertex $x \notin S$ and show that we can associate at least two detectors with $x$. For $v \in V(G)$, let $R_6(v) = N[v] \cup \{u \in V(G) : |N(u) \cap N(v)| = 2\}$. We will say that $x$ can only be associated with detector vertices within $R_6(x)$. Specifically, we will allow partial ownership of detectors, so a detector vertex, $v \in S$, contributes $\frac{1}{k}$, where $k = |R_6(v) \cap \overline{S}|$, toward the required total of two detectors.

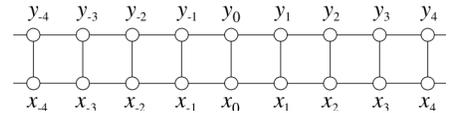

Figure 8: Ladder graph labeling

To begin, we will say that $x = x_0$, using the labeling convention shown in Figure 8. Suppose that $y_0 \notin S$; then $y_{-1}, y_1 \in S$ to 2-dominate $y_0$ and $x_{-1}, x_1 \in S$ to 2-dominate $x_0$. We see that $x_1$ and $y_1$ are not distinguished, so we need $x_2, y_2 \in S$; by symmetry, $x_{-2}, y_{-2} \in S$. Then $x$ receives a $\frac{1}{2}$ contribution from each of $\{x_{-1}, y_{-1}, x_1, y_1\}$, and we are done. Otherwise, we can assume $y_0 \in S$; suppose that $y_1 \notin S$. We require $x_1, x_2 \in S$ to 2-dominate $x_1$ and $y_{-1} \in S$ to 2-dominate $y_0$. Vertices $x_1$ and $x_2$ are not distinguished, so we need $x_3, y_2 \in S$, and by symmetry $x_{-1}, y_{-2} \in S$. Then $x$ receives at least $\frac{1}{2}$ from each of $\{x_1, y_0, x_{-1}, y_{-1}\}$, and we are done. Otherwise, we can assume $y_1 \in S$ and by symmetry $y_{-1} \in S$, in addition to $y_0 \in S$ that we showed previously. Vertex $x_0$ must be 2-dominated; without loss of generality let $x_1 \in S$. We see that $\{x_2, y_2\} \subseteq \overline{S}$ would cause $x_1$ and $y_2$ to not be distinguished, so $x_2$ and $y_2$ each contribute at least $\frac{1}{2}$ to $x$. If $x_{-1} \in S$, then $y_0$ contributes the final $\frac{1}{1}$ and we are done, so we assume that $x_{-1} \notin S$; then $y_0$ contributes $\frac{1}{2}$ and we need only another $\frac{1}{2}$ to have a total of two. We require $x_{-2} \in S$ to 2-dominate $x_{-1}$, and $y_{-2} \in S$ to distinguish $y_{-1}$ and $y_0$. Then $y_{-1}$ contributes the final $\frac{1}{2}$, completing the proof. □

**Theorem 9.** *For $j \geq 4$, we have RED:IC($P_2 \square P_j$) = RED:IC($P_2 \square C_j$) = $\lceil \frac{2}{3} n \rceil$.*

*Proof.* First, note that if $S_G$ and $S_H$ are RED:ICs on disjoint graphs $G$ and $H$, then $S_G \cup S_H$ is a RED:IC on $G + H + E_{GH}$ where $E_{GH}$ is any set of new edges between $G$ and $H$. Thus, if $S$ is a RED:IC set on $P_2 \square P_j$ with $|S| < \lceil \frac{2}{3} n \rceil$, then we can repeatedly create duplicates of $P_2 \square P_j$ (with duplicated detectors) and connect them end-to-end to produce a RED:IC set on $P_2 \square P_\infty$ with density strictly less than $\frac{2}{3} n$, which contradicts Theorem 8. Thus, RED:IC($P_2 \square P_j$) $\geq \lceil \frac{2}{3} n \rceil$. For the cylinder, $P_2 \square C_j$, we see that any maximal ladder subgraph is spanning and can be used as a tile to construct the infinite path, so we similarly find that RED:IC($P_2 \square C_j$) $\geq \lceil \frac{2}{3} n \rceil$. Figure 9 gives a set of RED:ICs on finite ladders and cylinders which achieve the lower bound of $\lceil \frac{2}{3} n \rceil$, completing the proof. □

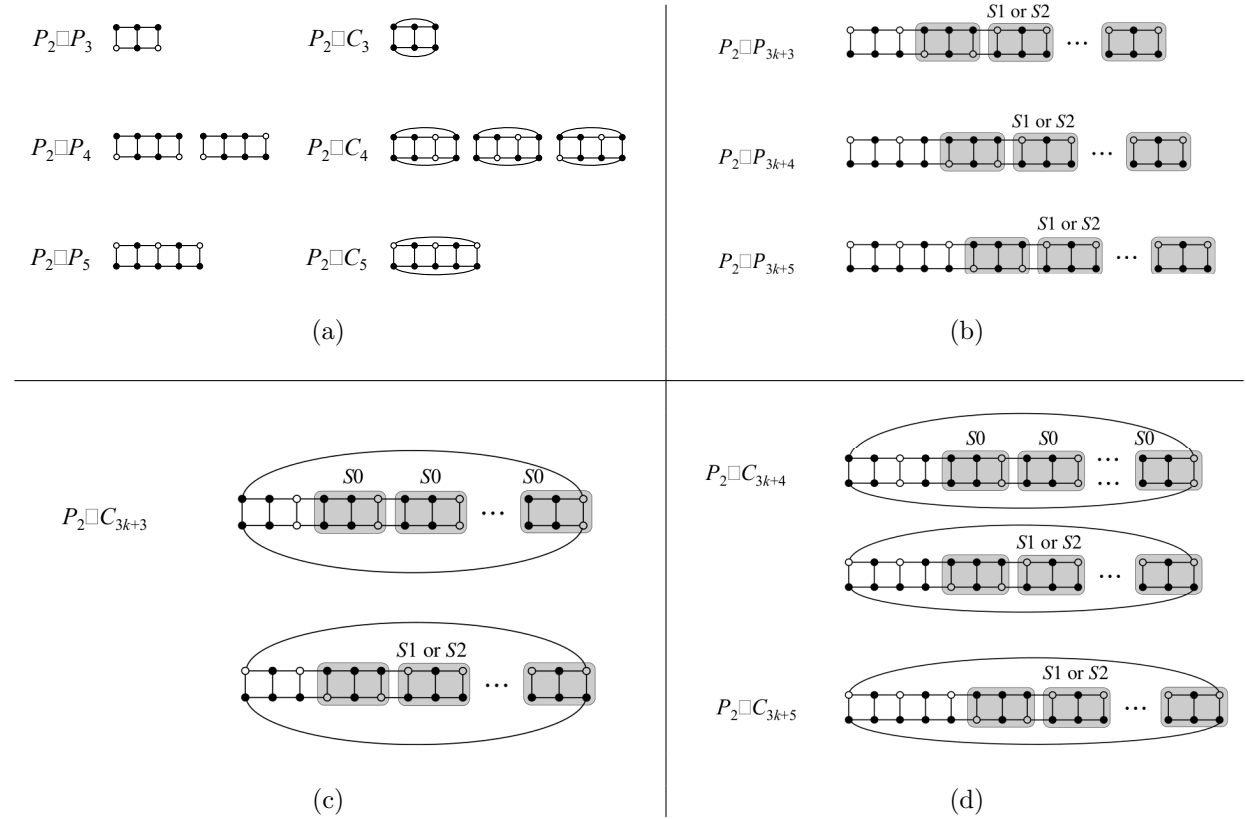

Figure 9: A set of optimal RED:ICs on finite ladders and cylinders. $S0$, $S1$, and $S2$ refer to $P_2 \square P_3$ tiles (disjoint subgraphs) which can be repeated. Base cases on $4 \leq n \leq 5$ vertices are given by (a).

**Theorem 10.** *For a finite torus $C_i \square C_j$, RED:IC($C_i \square C_j$) $\geq \lceil \frac{2}{5} n \rceil$.*

*Proof.* A torus can be drawn as $P_i \square P_j$, plus additional "wrapping" edges. We see that if $i \geq 5$, $j \geq 5$, and $i$ or $j$ is even, then any RED:IC set, $S$, on $C_i \square C_j$ can be tiled to produce a solution on the SQ. Because RED:IC%(SQ) $\geq \frac{2}{5}$, we know that $|S| \geq \lceil \frac{2}{5} n \rceil$.　□

### 4.3 Hypercubes

Let $Q_n = P_2^n$, where $G^n$ denotes repeated application of the $\square$ operator, be the hypercube in $n$ dimensions. If $S$ is a RED:IC on $Q_n$ for $n \geq 2$, then we can duplicate the vertices to produce a new RED:IC of size $2|S|$ on $Q_{n+1} = Q_n \square P_2$; thus, RED:IC%($Q_n$) is a non-increasing sequence in terms of $n$. We have found that RED:IC%($Q_5$) = $\frac{3}{8}$, which serves as an upper bound for the minimum density of RED:IC sets in larger hypercubes. Figure 10 shows an optimal RED:IC set for each of the hypercubes on $n \leq 5$ dimensions.

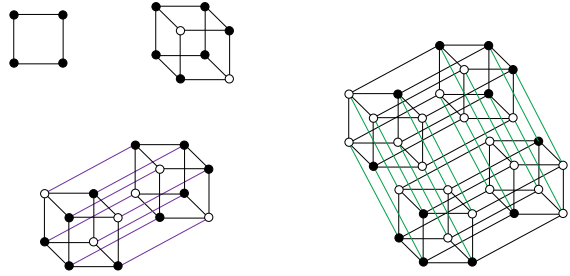

Figure 10: Optimal RED:IC for $Q_n$ with $n \leq 5$

### 4.4 Cubic Graphs

**Observation 7.** *RED:IC exists for all closed-twin-free cubic graphs.*

**Observation 8.** *On a cubic graph, RED:IC exists if and only if IC exists.*

#### Lower bound on RED:IC(G) for cubic

In the following subsections, we will frequently need to discuss sums of share values. To shorthand this process, we will let $\sigma_A$ for some sequence of single-character symbols, $A$, denote $\sum_{k \in A} \frac{1}{k}$. Thus, $\sigma_a = \frac{1}{a}$, $\sigma_{ab} = \frac{1}{a} + \frac{1}{b}$, and so on. For convenience of writing some expressions, we will also let the *domination number* of some $v \in V(G)$ be denoted by $dom(v)$; we define that $dom(v) = k$ if and only if $v$ is $k$-dominated.

**Theorem 11.** *If $G$ is a cubic graph, then RED:IC%($G$) $\geq \frac{4}{7}$.*

*Proof.* Let $x \in S$ be an arbitrary detector, and let $N(x) = \{a, b, c\}$. Among the three vertices $a$, $b$, and $c$, we have at most one edge, as otherwise we create closed-twins and a RED:IC would not exist. Suppose $ab \in E(G)$, or by symmetry $bc \in E(G)$ or $ac \in E(G)$. We know that $\exists z_1 \in N(a) - (N[b] \cup N[c])$ and $\exists z_2 \in N(b) - (N[a] \cup N[c])$, as otherwise we create closed-twins. In order to distinguish $x$ and $a$ we require $c, z_1 \in S$; and by symmetry to distinguish $x$ and $b$ we require $c, z_2 \in S$. If $a \in S$ or $b \in S$ then $sh(x) \leq \sigma_{3332} = \frac{1}{3} + \frac{1}{3} + \frac{1}{3} + \frac{1}{2} = \frac{3}{2} < \frac{7}{4}$, and we are done; otherwise $a, b \notin S$. To distinguish $x$ and $c$, we require $dom(c) = 4$, so $sh(x) \leq \sigma_{4222} = \frac{7}{4}$ and we are done. Otherwise, we can assume that there are no edges among $a$, $b$, and $c$. Suppose $dom(x) = 2$; let $a \in S$ and $b, c \notin S$. As seen in the previous case, to distinguish $x$ and $a$ we require $dom(a) = 4$, so $sh(x) \leq \sigma_{4222}$ and we are done. Similarly, if $dom(x) = 4$, then we are done, which leaves the last remaining case: $dom(x) = 3$. Let $a, b \in S$ and $c \notin S$. If $dom(a) \geq 3$ or $dom(b) \geq 3$ then $sh(x) \leq \sigma_{3322}$ and we are done; otherwise $dom(a) = dom(b) = 2$. We see that $x$ and $a$ are not distinguished, a contradiction. Thus, in any case $sh(x) \leq \frac{7}{4}$, giving a lower bound of RED:IC%($G$) $\geq \frac{4}{7}$ and completing the proof.　□

**RED:IC on the infinite 3-regular tree**

**Theorem 12.** *The infinite cubic tree, $T$, has RED:IC%$(T) = \frac{4}{7}$.*

*Proof.* Theorem 11 gives us a lower bound of RED:IC%$(T) \geq \frac{4}{7}$. The figure given in Figure 11 gives a RED:IC, $S$, on $T$. We see that every detector vertex, $x \in S$, has $sh(x) = \sigma_{4222} = \frac{7}{4}$, meaning the density of $S$ in $T$ is $\frac{4}{7}$, completing the proof. □

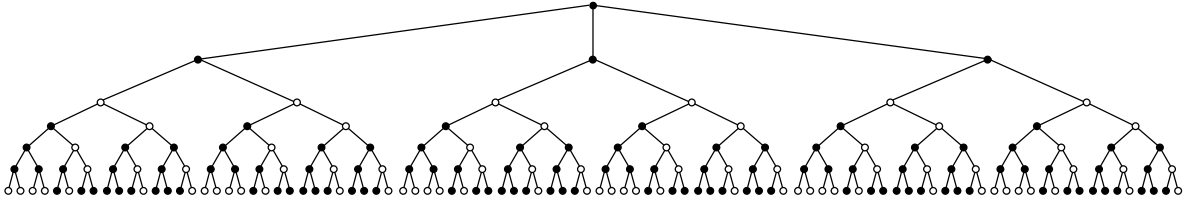

Figure 11: RED:IC%$(T) \leq \frac{4}{7}$

**RED:IC on the infinite hexagonal grid**

For the hexagonal grid, HEX, the tiling of the solution given in Figure 12 contains $\frac{2}{3}$ of the vertices as detectors; thus, we have RED:IC%$(HEX) \leq \frac{2}{3}$. The lower bound is from Theorem 11.

**Theorem 13.** *For the infinite hexagonal grid, HEX, $\frac{4}{7} \leq $ RED:IC%$(HEX) \leq \frac{2}{3}$.*

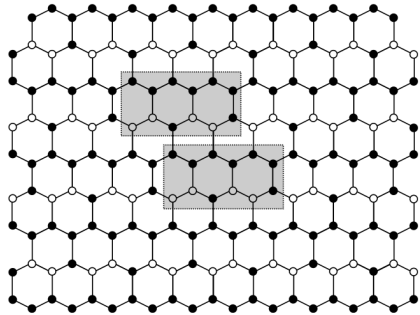

Figure 12: RED:IC%$(HEX) \leq \frac{2}{3}$

**Extremal cubic graphs with lower bound**

Let $G_{14}$ be the subgraph shown in Figure 13; $G_{14}$ contains four "loose" edges which may go to arbitrary vertices inside or outside of the subgraph (so long as the result is cubic). We see that each vertex in $G_{14}$ is at least 2-dominated, and it can be shown that each pair of vertices is 2-distinguished regardless of the specific incidence of the loose edges. For example, vertex pair $a, e$ is distinguished by $\{b, f\}$, vertex pair $b, e$ is distinguished by $\{b, c\}$, vertex pair $e, g$ is distinguished by $\{a, c\}$, and so on. Thus, we see that $G_{14}$ has a RED:IC of size 8. From Theorem 11, we know that a cubic graph must have RED:IC%$(G) \geq \frac{4}{7}$; thus, any cubic graph constructed using (only) copies of $G_{14}$ will have the minimum density of $\frac{4}{7}$. Copies of the $G_{14}$ subgraph can be connected in a ring to create an infinite family of cubic graphs which have the extremal value of RED:IC%$(G) = \frac{4}{7}$. This construction is shown in Figure 14.

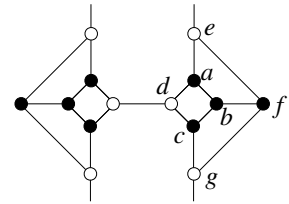

Figure 13: Cubic subgraph on 14 vertices requiring 8 detectors

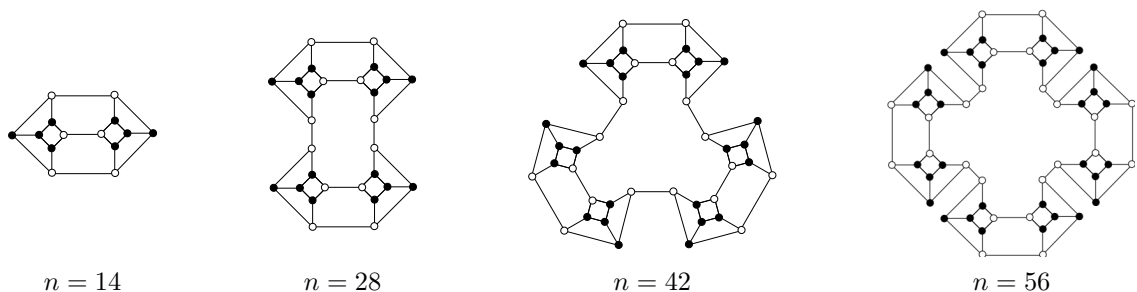

$n = 14$ $\qquad n = 28 \qquad n = 42 \qquad n = 56$

Figure 14: Infinite family of cubic graphs with RED:IC(G)=$\frac{4}{7}$

**Extremal cubic graphs with upper bound**

Let $G_6$ be the subgraph on 6 vertices from Figure 15; $G_6$ has two "loose" edges which extend out from $a$ and $d$ to any external vertex, so long as the entire graph is cubic. We see that vertices $b$ and $f$ can only be distinguished by having $\{c, e\} \subseteq S$, and by symmetry $\{b, f\} \subseteq S$. If $G$ is composed exclusively of disjoint copies of $G_6$ (allowing loose edges to overlap), then each vertex like $a$ or $d$ must be connected to another vertex like $a$ or $d$, as all other vertices already have degree three. We see that to distinguish $a$ and $b$, we require the vertex adjacent to $a$ by its loose edge to be a detector, so by symmetry all vertices like $a$ and $d$ must be detectors. Thus, all vertices in $G_6$ must be detectors. Copies of the $G_6$ subgraph can be connected in a ring to form an infinite family of cubic graphs with the extremal value RED:IC%$(G) = 1$, as shown in Figure 16.

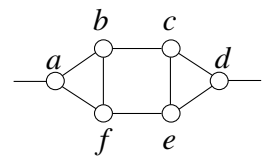

Figure 15: Cubic subgraph on 6 vertices requiring 6 detectors

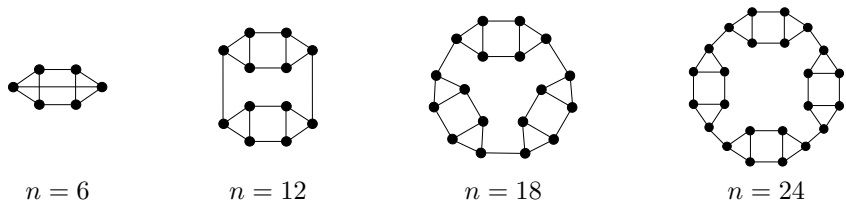

$n = 6 \qquad n = 12 \qquad n = 18 \qquad n = 24$

Figure 16: Infinite family of cubic graphs with RED:IC(G)=n

Table 2 gives a summary of results for the number of cubic graphs on up to 20 vertices which has a given value for RED:IC$(G)$.

| $n$ | 6 | 8 | 10 | 12 | 14 | 16 | 18 | 20 |
|---|---|---|---|---|---|---|---|---|
| cubic graphs | 2 | 5 | 19 | 85 | 509 | 4060 | 41301 | 510489 |
| with RED:IC | 2 | 4 | 14 | 63 | 386 | 3189 | 33586 | 427277 |
| lowest RED:IC($G$) | 6 | 6 | 6 | 8 | 8 | 10 | 11 | 12 |
| highest RED:IC($G$) | 6 | 6 | 8 | 12 | 12 | 14 | 18 | 18 |

Table 2: Results on RED:ICs for finite (connected) cubic graphs